\documentclass[11pt]{amsart}

\newcommand{\K}{\bold{K}}
\newcommand{\Q}{\bold{Q}}
\newcommand{\R}{\bold{R}}
\newcommand{\Z}{\bold{Z}}
\newcommand{\al}{\alpha}
\newcommand{\be}{\beta}

\newtheorem{lemma}{Lemma}
\newtheorem{theorem}{Theorem}

\theoremstyle{definition}
\newtheorem{conjecture}{Conjecture} 

\theoremstyle{remark}
\newtheorem{notation}{Notation} 
\newtheorem{note}{Note} 

\begin{document}

\title{Primitive divisors of Lucas and Lehmer sequences}
\author{Paul M Voutier}
\address{Probl\`{e}mes Diophantiens \\ 
	 Universit\'{e} P et M Curie (Paris VI) \\ 
	 Math\'{e}matiques, UFR 920 \\
	 Tour 45-46, 5\`{e}me Etage, B.P. 172 \\ 
	 4, place Jussieu \\ 
	 F-75252 PARIS Cedex 05, FRANCE} 
\email{voutier@@ccr.jussieu.fr}
\date{}
    
\maketitle

\begin{abstract}
Stewart reduced the problem of determining all Lucas and Lehmer 
sequences whose $n$-th element does not have a primitive divisor 
to solving certain Thue equations. Using the method of Tzanakis and 
de Weger for solving Thue equations, we determine such sequences 
for $n \leq 30$. Further computations lead us to conjecture that, 
for $n > 30$, the $n$-th element of such sequences always has a 
primitive divisor.  
\end{abstract}

\section{Introduction}

Let $\al$ and $\be$ be algebraic numbers such that $\al + \be$ 
and $\al \be$ are relatively prime non-zero rational integers 
and $\al / \be$ is not a root of unity. The sequence 
${ \left( u_{n} \right) }_{n=0}^{\infty}$ defined by 
$u_{n} = \left( {\al}^{n}-{\be}^{n} \right)/(\al - \be)$ 
for $n \geq 0$ is called a {\it Lucas sequence}. 

If, instead of supposing that $\al+\be \in \Z$, we only suppose that 
$(\al+\be)^{2}$ is a non-zero rational integer, still relatively 
prime to $\al \be$, then we define the {\it Lehmer sequence} 
${ \left( u_{n} \right) }_{n=0}^{\infty}$ associated to $\al$ and $\be$ by 
\begin{displaymath}
u_{n} = \begin{cases}
	   \displaystyle \frac{{\al}^{n}-{\be}^{n}}{\al - \be}  
	      & \text{if $n$ is odd}, \\
	   \displaystyle \frac{{\al}^{n}-{\be}^{n}}{{\al}^{2}-{\be}^{2}}
	      & \text{if $n$ is even}. 
	\end{cases}
\end{displaymath}

We say that a prime number $p$ is a {\it primitive divisor}
of a Lucas number $u_{n}$ if $p$ divides $u_{n}$ but does not
divide $(\al -\be)^{2}u_{2} \dotsm u_{n-1}$. Similarly, $p$ is a
primitive divisor of a Lehmer number $u_{n}$ if $p$ divides $u_{n}$
but not $(\al^{2} - \be^{2})^{2}u_{3} \dotsm u_{n-1}$.

There is another sequence, $(v_{n})_{n=0}^{\infty}$, 
associated to every Lucas and Lehmer sequence. However, 
$v_{n}=u_{2n}/u_{n}$ so it has a primitive divisor if and only if 
$u_{2n}$ does. Therefore, in what follows we need only consider 
the numbers $u_{n}$.

Drawing upon the ideas of Schinzel \cite{Sch} and refined techniques 
for determining lower bounds for linear forms in logarithms, Stewart 
\cite[p.80]{Ste1} showed, as a consequence of his Theorem 1, that if 
$n > C$ then $u_{n}$ has a primitive divisor, where $C=e^{452}2^{67}$ 
for Lucas sequences and $C=e^{452}4^{67}$ for Lehmer sequences.

Moreover, Stewart \cite[Theorem 2]{Ste1} also proved that for $n>6, 
n \neq 8,10,12$ there are at most finitely many Lucas and Lehmer 
sequences whose $n$-th element is without a primitive divisor.  
And, as Stewart states (see \cite[p.80]{Ste1}), for Lucas sequences 
the conditions on $n$ may be replaced by $n>4, n \neq 6$. He demonstrated 
this by reducing the problem of finding all Lucas and Lehmer sequences 
whose $n$-th element has no primitive divisor to solving finitely many 
Thue equations.  

Here we will use the method of Tzanakis and de Weger \cite{deW2,TW} 
to solve these Thue equations and thus enumerate all Lucas and Lehmer 
sequences whose $n$-th element has no primitive divisor when $n \leq 30$ 
satisfies the conditions in the previous paragraph.

Notice that if $(u_{n})$ is the sequence generated by $\al$ and 
$\be$ and $(u_{n}')$ is the sequence generated by $-\al$ and 
$-\be$, then $u_{n}=\pm u_{n}'$. We list only one of these pairs 
in Table~\ref{tab:Lucas} below, thus for each entry $(\al,\be)$ 
which generates a sequence whose $n$-th element is 
without a primitive divisor, the $n$-th element of the sequence  
generated by $-\al$ and $-\be$ also lacks a primitive divisor. 

Similarly, if $\al$ and $\be$ generate a Lehmer sequence $(u_{n})$ 
and $(u_{n}')$ is the Lehmer sequence generated by $i \al$ and 
$i \be$, then $u_{n} = \pm u_{n}'$. Again, we list only one of 
the four pairs in Table~\ref{tab:Lehmer} below, and so for each 
entry $(\al,\be)$ generating a sequence whose $n$-th element is 
without a primitive divisor, the $n$-th element of the sequences 
generated by $i^{k} \al$ and $i^{k} \be$, for $k=1,2,3$, also lack 
primitive divisors. 

\begin{theorem}
{\rm (i)} For $4 < n \leq 30$, $n \neq 6$, Table~{\rm \ref{tab:Lucas}} 
gives a complete list, up to the sign of $\al$ and $\be$, of all 
Lucas sequences whose $n$-th element has no primitive divisor.  

{\rm (ii)} For $6 < n \leq 30$, $n \neq 8,10,12$, Table~{\rm \ref{tab:Lehmer}}
gives a complete list, up to multiplication 
of $\al$ and $\be$ by a fourth-root of unity, of all Lehmer sequences 
whose $n$-th element has no primitive divisor. 
\end{theorem}

\begin{table}
\caption{}
\label{tab:Lucas}
\begin{center}
\begin{tabular}{||c|c|c|c|c||}                         \hline
$n$  & \multicolumn{4}{c||}{$\al, \be$}             \\ \hline
$5$  & $\displaystyle \frac{1 \pm \sqrt{5}}{2}$        
     & $\displaystyle \frac{1 \pm \sqrt{-7}}{2}$       
     & $ 1 \pm \sqrt{-10} $                            
     & $\displaystyle \frac{1 \pm \sqrt{-11}}{2}$   \\ \cline{2-5}
     & $\displaystyle \frac{1 \pm \sqrt{-15}}{2}$      
     & $ 6 \pm \sqrt{-19} $                         
     & $ 6 \pm \sqrt{-341}$                          
     &                                              \\ \hline
$7$  & $\displaystyle \frac{1 \pm \sqrt{-7}}{2}$       
     & $\displaystyle \frac{1 \pm \sqrt{-19}}{2}$      
     & &                                            \\ \hline 
$8$  & $ 1 \pm \sqrt{-6} $                             
     & $\displaystyle \frac{1 \pm \sqrt{-7}}{2}$       
     & &                                            \\ \hline 
$10$ & $ 1 \pm \sqrt{-2} $                             
     & $\displaystyle \frac{5 \pm \sqrt{-3}}{2}$       
     & $\displaystyle \frac{5 \pm \sqrt{-47}}{2}$      
     &                                              \\ \hline
$12$ & $\displaystyle \frac{1 \pm \sqrt{5}}{2}$     
     & $\displaystyle \frac{1 \pm \sqrt{-7}}{2}$       
     & $\displaystyle \frac{1 \pm \sqrt{-11}}{2}$   
     &                                              \\ \cline{2-5}   
     & $ 1 \pm \sqrt{-14} $                            
     & $\displaystyle \frac{1 \pm \sqrt{-15}}{2}$      
     & $\displaystyle \frac{1 \pm \sqrt{-19}}{2}$   
     &                                              \\ \hline
$13$ & $\displaystyle \frac{1 \pm \sqrt{-7}}{2}$       
     & & &                                          \\ \hline
$18$ & $\displaystyle \frac{1 \pm \sqrt{-7}}{2}$       
     & & &                                          \\ \hline
$30$ & $\displaystyle \frac{1 \pm \sqrt{-7}}{2}$      
     & & &                                          \\ \hline
\end{tabular}
\end{center}
\end{table}

\begin{table}
\caption{}
\label{tab:Lehmer}
\begin{center}
\begin{tabular}{||c|c|c|c||}                                  \hline
$n$  & \multicolumn{3}{c||}{$\al, \be$}                    \\ \hline
$7$  & $\displaystyle \frac{1 \pm \sqrt{-7}}{2}$  
     & $\displaystyle \frac{1 \pm \sqrt{-19}}{2}$  
     & $\displaystyle \frac{\sqrt{3} \pm \sqrt{-5}}{2}$    \\ \cline{2-4}
     & $\displaystyle \frac{\sqrt{5} \pm \sqrt{-7}}{2}$    
     & $\displaystyle \frac{\sqrt{13} \pm \sqrt{-3}}{2}$   
     & $\displaystyle \frac{\sqrt{14} \pm \sqrt{-22}}{2}$  \\ \hline
$9$  & $\displaystyle \frac{\sqrt{5} \pm \sqrt{-3}}{2}$           
     & $\displaystyle \frac{\sqrt{7} \pm \sqrt{-1}}{2}$           
     & $\displaystyle \frac{\sqrt{7} \pm \sqrt{-5}}{2}$    \\ \hline
$13$ & $\displaystyle \frac{1 \pm \sqrt{-7}}{2}$ 
     & &                                                   \\ \hline
$14$ & $\displaystyle \frac{\sqrt{3} \pm \sqrt{-13}}{2}$   
     & $\displaystyle \frac{\sqrt{5} \pm \sqrt{-3}}{2}$   
     & $\displaystyle \frac{\sqrt{7} \pm \sqrt{-1}}{2}$    \\ \cline{2-4}
     & $\displaystyle \frac{\sqrt{7} \pm \sqrt{-5}}{2}$    
     & $\displaystyle \frac{\sqrt{19} \pm \sqrt{-1}}{2}$    
     & $\displaystyle \frac{\sqrt{22} \pm \sqrt{-14}}{2}$  \\ \hline
$15$ & $\displaystyle \frac{\sqrt{7} \pm \sqrt{-1}}{2}$  
     & $\displaystyle \frac{\sqrt{10} \pm \sqrt{-2}}{2}$  
     &                                                     \\ \hline
$18$ & $\displaystyle \frac{1 \pm \sqrt{-7}}{2}$   
     & $\displaystyle \frac{\sqrt{3} \pm \sqrt{-5}}{2}$    
     & $\displaystyle \frac{\sqrt{5} \pm \sqrt{-7}}{2}$    \\ \hline
$24$ & $\displaystyle \frac{\sqrt{3} \pm \sqrt{-5}}{2}$  
     & $\displaystyle \frac{\sqrt{5} \pm \sqrt{-3}}{2}$  
     &                                                     \\ \hline
$26$ & $\displaystyle \frac{\sqrt{7} \pm \sqrt{-1}}{2}$           
     & &                                                   \\ \hline
$30$ & $\displaystyle \frac{1 \pm \sqrt{-7}}{2}$    
     & $\displaystyle \frac{\sqrt{2} \pm \sqrt{-10}}{2}$   
     &                                                     \\ \hline
\end{tabular}
\end{center}
\end{table}

Using continued-fractions and Lemma~\ref{lem:start}(i) below, we 
may quickly search for small solutions of Thue equations. By this 
method, we have determined that for $31 \leq n \leq 250$ there are 
no solutions $(x,y)$ of the appropriate Thue equations, $F_{n}(X,Y)=m$ 
with $\max (|x|,|y|) < 10^{6}$ which give rise to Lucas or Lehmer 
sequences whose $n$-th element is without a primitive divisor. 
Notice that $10^{6}$ is quite a bit larger than the entries $X_{4}$ 
and $Y_{4}$ in Tables 5--7 below, which give the maximum of $|x|$ 
and $|y|$ for all solutions $(x,y)$ of each completely solved Thue 
equation. Moreover, Birkhoff and Vandiver \cite{BV} have shown that 
there are no Lucas sequences generated by $\al, \be \in \Z$ whose 
$n$-th element does not have a primitive divisor for $n > 6$ and 
Carmichael \cite{Car} proved the same result for $\al, \be \in \R$ 
with $n > 12$. In the late 1950's, Ward \cite{Ward} and Durst 
\cite{Durst} extended Carmichael's result to Lehmer sequences. 
Thus it seems reasonable to make the following conjecture. 

\begin{conjecture}
For $n > 30$, the $n$-th element of a Lucas or Lehmer sequence 
always has a primitive divisor.
\end{conjecture}

In the next section, the lemmas necessary to establish the 
connection between enumerating Lucas and Lehmer sequences whose 
$n$-th element has no primitive divisor and solving Thue equations 
are given as well as some results from algebraic number theory 
which are necessary to solve these equations. In \S~3 we 
consider the values $n=5,8,10$ and $12$ for which there are 
infinitely many Lehmer sequences without a primitive divisor but 
only finitely many such Lucas sequences. Then, in \S~4, we give 
a description of the algorithm of Tzanakis and de Weger for solving 
Thue equations. In \S~5, we describe linear dependence relations 
over $\Z$ between certain numbers which arise in our applications of 
this algorithm. Finally, tables giving details of the computations 
for each value of $n$ are provided.  

This paper originated from my Master's thesis conducted under the 
supervision of Dr. C. L. Stewart at the University of Waterloo. 
Dr. Stewart deserves my deepest thanks for his encouragement,
patience and knowledge during my studies under him. This work was
completed with the help of a Graduate School Dean's Small Grant
Award from the University of Colorado, Boulder as well as the 
support provided by my Ph.D. advisor Dr. W. M. Schmidt. 

\section{Some Preliminary Lemmas}

Let $\Phi_{n}(X,Y)$ be the homogeneous cyclotomic polynomial of 
order $n$ and $\phi_{n}(X) = \Phi_{n}(X,1)$. These polynomials 
are linked to Lucas and Lehmer sequences by the formula
\begin{equation}
\label{eq:connect}
{\al}^{n}-{\be}^{n} = \prod_{d|n}\Phi_{d}(\al, \be). 
\end{equation}

\begin{notation}
For $n > 1$, we let $P(n)$ denote the largest prime divisor of $n$.
\end{notation}

\begin{lemma}
\label{lem:tophi}
Let $n>4$ and $n \neq 6,12$. Then $u_{n}$ has a primitive divisor 
if and only if $\Phi_{n}(\al, \be) \neq \pm 1, \pm P(n/(n,3))$. 
$u_{12}$ has a primitive divisor if and only if 
$\Phi_{12}(\al, \be) \neq \pm 1, \pm 2, \pm 3, \pm 6$. 
\end{lemma}

\begin{proof}
This follows immediately from Lemmas 6 and 7 of Stewart 
\cite{Ste2}, using (\ref{eq:connect}). 
\end{proof}

It is by means of this lemma that we obtain Thue equations. For 
\begin{displaymath}
\Phi_{n}(\al, \be) = \prod_{\substack{
			       j=1 \\ 
			       (j,n)=1 
			    }  
			   }^{n}
		     (\al -\zeta_{n}^{j} \be )
= \prod_{\substack{
	    j=1 \\ 
	    (j,n)=1
	 }   
	}^{n/2} 
  \left( {\al}^{2}+{\be}^{2}
  - \left( \zeta_{n}^{j}+\zeta_{n}^{-j} \right) \al \be
\right).  
\end{displaymath}
Let $x={\al}^{2}+{\be}^{2}$ and $y=\al \be$. Then  
$\Phi_{n}(\al,\be)=F_{n}(x,y)$ where $F_{n}(x,y)$ is 
a binary form of total degree $\varphi (n)/2$ in $x$ 
and $y$ and $\varphi (n)$ is the Euler phi-function. 
Moreover, since $\zeta_{n}^{j}+\zeta_{n}^{-j}=2\cos(2\pi j/n)$ 
is an algebraic integer, $F_{n}(x,y)$ has rational integer 
coefficients. By Lemma~\ref{lem:tophi}, we have four Thue 
equations associated with each $n > 4$, $n \neq 6,12$ to 
solve. Recall that we assumed $(\al +\be)^{2}$ and $\al \be$ 
are integers so $x$ and $y$ are integers. We can find the 
values of $\al$ and $\be$ associated to a given solution 
$(x,y)$ of these Thue equations as the roots of the polynomial 
$X^{2} - \sqrt{x+2y} \, X + y$. 

We start with some properties of these binary forms.  

\begin{lemma}
\label{lem:sqft}
{\rm (i)} If $t$ is an odd integer then $F_{2t}(X,Y)=F_{t}(X,-Y)$. 

{\rm (ii)} Let $n=p_{1}^{r_{1}} \dotsm p_{k}^{r_{k}}$ and 
$m=p_{1}^{s_{1}} \dotsm p_{k}^{s_{k}}$, where $p_{1}, \dots, p_{k}$ 
are distinct primes and the $r_{i}$'s and $s_{i}$'s are positive 
integers with $1 \leq s_{i} \leq r_{i}$. Then $F_{n}(X,Y)=F_{m}(X',Y')$, 
where $Y'=Y^{n/m}$ and $X'$ can be written as a binary form of degree 
$n/m$ in $X$ and $Y$ with integer coefficients. 
\end{lemma}

\begin{proof}
These two statements follow easily from the analogous 
statements which hold for the cyclotomic polynomials:  
\begin{displaymath}
\phi_{2t}(X) = \phi_{t}(-X) \hspace{ 5mm} \text{ and }  
\hspace{5 mm} \phi_{n}(X) = \phi_{m} \left( X^{n/m} \right). 
\end{displaymath}

The first result is part (iv) of Prop.\ 5.16 from Chapter 2 
of Karpilovsky's book \cite{Karp}, while the second is a 
slight generalization, whose proof is essentially identical, 
of part (vi) of this same proposition.    
\end{proof}

There are two other results we need in order to implement the 
algorithm of Tzanakis and de Weger. All but one of the Thue 
equations we consider here split into linear factors in the 
field $\Q(\cos(2\pi/n))$. We need a factorization of the ideal 
$(P(n/(3,n)))$ in these fields as well as a system of fundamental 
units for the ring of integers of these fields.

\begin{lemma}
\label{lem:factor}
{\rm (i)} If n=$p^{k}$ is an odd prime power satisfying 
	  $7 \leq p^{k} \leq 29$, then 
\begin{displaymath}
(P(n/(3,n))) = (p) = (2-2\cos(2\pi/p^{k}))^{\varphi (p^{k})/2}. 
\end{displaymath}

{\rm (ii)} For $n = 15,21$ and $24$,
	   $(P(n/(3,n))) = (1+2\cos(2\pi/n))^{\varphi (n)/2}$. 
     
{\rm (iii)} For $n = 16$ and $20$, 
	    $(P(n/(3,n)))=(2\cos(2\pi/n))^{\varphi (n)/2}$.  
\end{lemma}

\begin{proof}
These factorizations are determined by using a theorem of Dedekind, 
see Proposition 2.14 of Washington \cite{Wash}. 
\end{proof}

\begin{note}
Notice that as a consequence of this lemma, any algebraic integer 
in $\Q(\cos(2\pi/n))$, where $n$ is as stated, with norm equal 
to $P(n/(3,n))$ must be an associate of the generator of the 
ideal given in these factorizations. 
\end{note}

\begin{lemma}
\label{lem:units}
{\rm (i)} If $n=p^{k}$, where $p$ is a prime and $\varphi (n) \leq 66$ 
	  then 
\begin{displaymath}
\left\{ \frac{\sin(a\pi/n)}{\sin(\pi/n)} : 1 < a < n/2, (a,n)=1 \right\} 
\end{displaymath}
is a system of fundamental units for $\Q( \cos(2\pi/n))$. 

In {\rm (ii)--(v)}, put $r=2\cos(2\pi/n)$. 

{\rm (ii)} If $n=15$ then $\{ r,r-1,r^{2}-3 \}$ is a system of 
fundamental units. 

{\rm (iii)} If $n=20$ then $\{ r-1,r-2,r^{2}-2 \}$ is a system of 
fundamental units.

{\rm (iv)} If $n=21$ then $\{ r,r-1,r^{2}+r-1,r^{2}-2,r^{2}-3 \}$ is a 
system of fundamental units.

{\rm (v)} If $n=24$ then $\{ r,2r-1,r^{2}-r-1 \}$ is a system of 
fundamental units.
\end{lemma}

\begin{proof}
(i) follows from Theorem 8.2 of Washington \cite{Wash} 
which states that the index of the group generated by these units 
in the full unit group is the class number and Theorem 1 of van 
der Linden \cite{vdL} which states that the class number of these 
fields is 1. 

(ii)--(v) The systems of units given in these cases were found 
by the methods of Pohst \& Zassenhaus \cite{PZ}. 
\end{proof}

\section{$n=5,8,10$ and $12$}

For $n = 5,8,10,12$, $F_{n}(X,Y)=m$ is of total degree two and 
reducible to a Pell equation. Hence there will be infinitely 
many solutions (if any) to the Thue equations which arise in these 
cases. However, for Lucas sequences not only is $(\al +\be)^{2}$ 
an integer but so is $\al +\be$, $x+2y$ must therefore be a 
perfect square for any solution $(x,y)$ of $F_{n}(X,Y)=m$. 
Letting $Z^{2}=X+2Y$ and substituting this expression for $X$ 
in $F_{n}(X,Y)=m$ we get an equation which can be transformed 
into one of the form $aX^{2}-bY^{4}=c$ where $a,b,c$ are pairwise 
relatively prime integers with $a$ and $b$ positive. 

We start with a lemma stating two results which will be used 
below and are likely to be difficult for the reader to find. 
First, a little notation. For a positive odd square-free number 
$A$, let $(a,b)$ be the least positive integer solution of 
$AX^{2}-Y^{2}=2$ and put $\K_{1}=\Q(\sqrt{A}), 
\K_{2}=\Q(\sqrt{b+a\sqrt{A}})$ and $\K_{3}=\Q(\sqrt{-b+a\sqrt{A}})$. 
We also define 
${\mathcal U}_{2} = \{ \al \in \K_{2}: {\mathcal N}_{K_{2}/K_{1}} (\al) = 1 \}
= \{ \pm \epsilon_{1}^{n}: n \in \Z \}$ and 
${\mathcal U}_{3} = \{ \al \in \K_{3}: {\mathcal N}_{K_{3}/K_{1}} (\al) = 1 \}
= \{ \pm \epsilon_{2}^{n}: n \in \Z \}$. For $\alpha \in \K_{2}$  
or $\K_{3}$, we let $\alpha'$ be its real conjugate.

\begin{lemma}
\label{lem:ljung}
{\rm (i)} Let $A$ be a square-free odd positive integer. For $x$ and $y$, 
positive integers satisfying $Ax^{2}-y^{4}=2$, we put 
$\vartheta = 2(Ax^{2}+y^{4}+2xy^{2}\sqrt{A})$. Then either 
\begin{displaymath}
\vartheta = \epsilon_{1} + \epsilon_{1}' + 2 \text{ or }  
\vartheta = \epsilon_{2} + \epsilon_{2}' + 2. 
\end{displaymath}

{\rm (ii)} Let $c$ and $D$ be positive squarefree integers with 
$(c,2)=(c,3)=(D,3)=1$. The three equations 
$8c^{2}X^{4} \pm 4cX^{2} + 1 = DY^{2}, c^{2}X^{4} + 1 = 2DY^{2}$ 
have between them at most one solution in positive integers. 
\end{lemma}

\begin{proof}
These results, both due to Ljunggren, are Satz 2 of 
\cite{Ljung1} and Satz VIII of \cite{Ljung2} respectively. 
\end{proof}

The author apologizes for the rather ad hoc form of this section, 
but there seems to be no way around this. There is, presently, no 
elegant and unified theory for solving these equations. One could 
use the Thue equation approach, used in the case $n=12,k=-2$, 
however there are apparently formidable computational difficulties 
for some of the Thue equations which would arise. 

\subsection{$n=5$}

In the case of $n=5$, we have 
$F_{5}(X,Y) = X^{2}+XY-Y^{2} = \pm 1, \pm 5$. Letting $X+2Y=Z^{2}$ 
where $Z$ is an integer, we obtain the equation 
$G_{5}(Y,Z) = Z^{4}-3YZ^{2}+Y^{2}=k=\pm 1, \pm 5$. Solving 
$G_{5}(Y,Z)-k=0$ as a quadratic in $Y$, we find that  
$2Y=3Z^{2} \pm \sqrt{5Z^{4}+4k}$, so that 
\begin{displaymath}
5Z^{4}+4k=W^{2}.       
\end{displaymath}

It is this last equation we shall use. 

In the cases of $k=\pm 1$, Cohn \cite[Theorem 7]{Cohn} has shown  
that the complete solution of $5Z^{4} \pm 4 = W^{2}$ in non-negative
integers is $(W,Z) = \{(1,1),(2,0),(3,1),(322,12)\}$. From 
this, we can show that the complete solution of $F_{5}(X,Y)=\pm 1, 
X+2Y=Z^{2}$ is $(X,Y) \in \{ (-610,377),(-5,3),(-3,2),(-2,1),
(-1,1),(1,0),(2,-1),(34,55) \}$. 

If $k=\pm 5$, we may substitute $W=5V$ and get the equation 
$Z^{4}+4=5V^{2}$. Here Cohn \cite[Theorem 13]{Cohn} has shown 
that the complete solution in non-negative integers is 
$(V,Z) \in \{(1,1),(2,2)\}$. Thus we find that the complete 
solution of $F_{5}(X,Y)=5, X+2Y=Z^{2}$ is 
$(X,Y) \in \{ (-18,11),(-7,4),(2,1),(3,-1) \}$. 

The Lucas sequences which arise from these solutions are given 
in Table~\ref{tab:Lucas}. Notice that some of these solutions 
do not give rise to Lucas sequences since for the corresponding 
values of $\al$ and $\be$, we have $\al \be = 0$ or $\al / \be$ 
is a root of unity. 

\subsection{$n=8$}

Here we obtain $2Z^{4}+2k=W^{2}$ where $k = \pm 1, \pm 2$ 
and $Y=(2Z^{2} \pm \sqrt{2Z^{4}+2k})/2$. As $W=2V$, we get 
the equation 
\begin{displaymath}
Z^{4}+k=2V^{2}.  
\end{displaymath}

For $k=1$, our equation is the third equation in 
Lemma~\ref{lem:ljung}(ii) with $c=D=1$. Notice that $(V,Z)=(1,1)$ 
is a solution in positive integers and, hence by this lemma, the 
only such solution. This implies that the complete solution of 
$F_{8}(X,Y)=1$ where $X+2Y$ is a perfect square is 
$(X,Y) \in \{ (-3,2),(1,0) \}$.     

If $k=-1$, we can factor $Z^{4}-1$, getting 
$(Z^{2}-1)(Z^{2}+1)=2V^{2}$. Notice that $Z$ must be odd,  
so that $(Z^{2}-1,Z^{2}+1)=2$. Therefore, one of $Z^{2} \pm 1$ 
is a square and the other twice a square. For  $Z \neq 0$, 
$Z^{2}+1$ is not a square and hence $Z^{2}-1$ must be a square. 
So $Z=0, \pm 1$ and we find that the complete solution of 
$F_{8}(X,Y)=-1, X+2Y=Z^{2}$ is $(X,Y) \in \{ (-1,1) \}$. 

If $k=2$, we see that $Z=2U$ and obtain the equation 
$8U^{4}=V^{2}-1=(V+1)(V-1)$. Notice that $V$ must be odd and 
so $(V-1,V+1)=2$. Thus either $V-1=2H^{4}$ and $V+1=4K^{4}$ 
or $V-1=4K^{4}$ and $V+1=2H^{4}$, from which we deduce that 
$H^{4}-2K^{4}=\pm 1$. Delone and Faddeev \cite[Theorem 3, p.374]{DF} 
have shown that the complete solution in non-negative integers of 
these equations is $(H,K) \in \{(1,1),(1,0)\}$. So we determine that 
the complete solution of $F_{8}(X,Y)=2$ with $X+2Y$ a perfect 
square is $(X,Y) \in \{ (-10,7),(-2,1),(2,-1),(2,1) \}$. 

If $k=-2$, the relation $Z=2U$ gives us the equation 
$V^{2}-8U^{4}=-1$ which has no solution, since $-1$ is 
not a square mod $8$. 

Computing the values of $\al$ and $\be$ which correspond to 
these solutions, we find that Table~\ref{tab:Lucas} is 
complete for $n=8$. 

\subsection{$n=10$}

Using the same argument as in the previous section, we find that 
$10Y=5Z^{2} \pm \sqrt{5Z^{4}+20k}$. If $W^{2}=5Z^{4}+20k$ then 
$W=5V$ so we consider the equations 
\begin{displaymath}
5V^{2}=Z^{4}+4k
\end{displaymath}
where again $k = \pm 1, \pm 5$. 

Notice that for $k = \pm 1$, we have the same equations as for 
$n=5$ and $k = \pm 5$. So we find the same values of $Z$. Here 
$(X,Y) \in \{ (-2,3),(-1,1),(1,0),(2,1) \}$ is the complete 
solution of $F_{10}(X,Y)=\pm 1, X+2Y=Z^{2}$. 

For $k=\pm 5$, we have $Z=5U$ which leads to the equation 
$V^{2}=125U^{4} \pm 4$. Using the theory of Pell equations,  
any solution must be the square root of five times a Fibonacci 
number. Robbins \cite[Theorem 3]{Rob} has shown that the only 
squares of the form five times a Fibonacci number are 0 and 25. 
Thus the complete solution of $F_{10}(X,Y)= \pm 5, X+2Y=Z^{2}$ 
is $(X,Y) \in \{ (-11,18),(-2,1),(2,-1),(11,7) \}$. 

Again, the Lucas sequences which arise from these solutions 
are given in Table~\ref{tab:Lucas}. 

\subsection{$n=12$}

We have $F_{12}(X,Y) = X^{2}-2Y^{2} = k = \pm 1,\pm 2,\pm 3,\pm 6$ 
and $Y = 2Z^{2} \pm \sqrt{3Z^{4}+k}$, so that 
\begin{displaymath}
3Z^{4}+k=W^{2}.   
\end{displaymath}

First we note that the equations for $k=-1,2,3,-6$ have no 
solutions by considering them mod $3$. 

Consider next the case of $k=1$. Here there are precisely two 
solutions of $3Z^{4}+1=W^{2}$ in positive integers, namely 
$(W,Z)=(1,2)$ or $(2,7)$ (see Ljunggren \cite{Ljung3}). Thus 
the complete solution of $F_{12}(X,Y)=1$ with $X+2Y=Z^{2}$ is 
$(X,Y) \in \{ (-26,15),(-7,4),(-2,1),(1,0),(2,-1),(2,1) \}$. 

If $k=-2$, we convert the equation into a Thue equation. We 
factor $W^{2}+2=3Z^{4}$ as $(W-\sqrt{-2})(W+\sqrt{-2}) = 3Z^{4}$. 
Thus we can write $W-\sqrt{-2} = AU^{4}$ and $W+\sqrt{-2} = BV^{4}$, 
with $AB=3L^{4}$ where $A,B,L,U$ and $V$ are algebraic integers in 
$\Q(\sqrt{2})$, with $A$ and $B$ algebraic conjugates, $U=S-T\sqrt{-2}$ 
and $V=S+T\sqrt{-2}$. Using the facts that $3=(1+\sqrt{-2})(1-\sqrt{-2})$, 
the only units in $\Q(\sqrt{-2})$ are $\pm 1$ and 
$(W-\sqrt{-2},W+\sqrt{-2})$ divides $(2\sqrt{-2})$, we find that $S$ 
and $T$ give rise to a solution of $W^{2}+2=3Z^{4}$ if and only if 
\begin{displaymath}
F(S,T) = S^{4}-4S^{3}T-12S^{2}T^{2}+8ST^{3}+4T^{4} = 1. 
\end{displaymath}

We will use the method of Tzanakis and de Weger to solve this 
Thue equation, finding that $(S,T) = ( \pm 1, 0)$ are the only 
solutions. To apply this method, we need a system of fundamental 
units for the ring of integers of $\Q(\al)$ where $F(\al,1)=0$. 
We find, by the methods of Pohst \& Zassenhaus \cite{PZ}, that 
$\{ \epsilon_{1}=(\al^{3}-4\al^{2}-10\al+12)/4,
\epsilon_{2}=(2\al^{3}-9\al^{2}-20\al+26)/4,
\epsilon_{3}=(2\al^{3}-9\al^{2}-24\al+26)/4 \}$ 
is such a system. We give some of the details of the computations 
used to solve this Thue equation in Table 5 of Section 6. From these 
solutions of the Thue equation, we find that $(W,Z) = ( \pm 1, \pm 1)$ 
are the only solutions of $W^{2}+2=3Z^{4}$ and thus the complete 
solution of $F_{12}(X,Y)=-2$ with $X+2Y=Z^{2}$ is 
$(X,Y) \in \{(-5,3),(-1,1)\}$. 

If $k=-3$, then $W=3V$ and we obtain the equation $Z^{4}-3V^{2}=1$. 
Ljunngren \cite{Ljung3} has shown that $(V,Z)=(0,\pm 1)$ are the 
only integer solutions. Therefore the complete solution of 
$F_{12}(X,Y)=-3$ with $X+2Y$ a perfect square is $(X,Y) = (-3,2)$. 

For $k=6$, $W=3V$ which leads to the equation $3V^{2}-Z^{4}=2$, 
the special case $A=3$ of the equation in Lemma~\ref{lem:ljung}(i). 
The methods of Pohst \& Zassenhaus \cite{PZ} show that 
$\epsilon_{1}+\epsilon_{1}'+2 = 24+12\sqrt{3}$, which gives  
arise to no solutions, and 
$\epsilon_{2}+\epsilon_{2}'+2 = 8+4\sqrt{3}$, which gives rise 
to the solutions $(V,Z)=(\pm 1, \pm 1)$. Thus the complete 
solution of $F_{12}(X,Y)=6$ with $X+2Y$ a perfect square is 
$(X,Y) \in \{ (-9,5),(3,-1) \}$. 

Again, the Lucas sequences which arise from these solutions 
are given in Table~\ref{tab:Lucas}. 

\section{The Algorithm}

We now describe the algorithm which is used to solve these 
Thue equations. We will follow the notation and numbering 
of constants from the paper of Tzanakis and de Weger \cite{TW}, 
except that we use $d$ where they use $n$ and label our 
linear forms more explicitly. The reader will find in their 
paper, proofs of the lemmas below as well as explicit formulas 
for the constants mentioned below.  

Suppose we wish to solve $F(X,Y)=m$ where 
\begin{displaymath}
F(X,Y) = f_{d}X^{d} + f_{d-1}X^{d-1}Y + \dots + f_{0}Y^{d} \in \Z [X,Y]
\end{displaymath}
is an irreducible polynomial of total degree $d \geq 3$ and 
$m$ is a non-zero integer. Let $\xi = \xi^{(1)}, \dots, 
\xi^{(d)}$ be the roots of $g(X)=F(X,1)$ which we shall 
assume are all real (this is not necessary but the polynomials 
which we shall consider have only real roots and it makes the 
notation easier) and let $\K = \Q \left( \xi \right)$.  
Notice that  $[\K : \Q] = d$. 

Let $\{ \epsilon_{1}, \dots , \epsilon_{d-1} \}$ be a system 
of fundamental units of ${\mathcal O}_{K}$, the ring of integers 
in $\K$. Let ${\mathcal D}$ be a full module in $\K$ containing 
$1$ and $\xi$. We can partition the set of elements $\mu$ of 
${\mathcal D}$ with $f_{d}{\mathcal N}_{K/Q}(\mu)=m$ into finitely 
many equivalence classes under the relation of being associates 
(see the corollary to Theorem 5 on p. 90 of \cite{BS}) and 
let ${\mathcal M}$ be a complete set of representatives of these 
equivalence classes. Notice that if $|f_{d}|=|m|=1$, then 
${\mathcal M} = \{ 1 \}$. By Lemma~\ref{lem:factor}, for our 
applications, when $m \neq \pm 1$, ${\mathcal M}$ contains only 
the generator of the ideal specified in there. 

\subsection{An Initial Upper Bound via Linear Forms in Logarithms} 

For $(x,y) \in {\Z}^{2}$ satisfying $F(x,y)=m$, we put 
$\be = x - \xi y$. For each $i$ and any choice of $j,k$ 
satisfying $i \neq j \neq k \neq i$, we let  
\begin{displaymath}
\Lambda(i,j,k,\mu) = \log \left| 
\frac{\xi^{(i)}-\xi^{(j)}}{\xi^{(i)}-\xi^{(k)}}
\frac{{\be}^{(k)}}{{\be}^{(j)}} \right|, 
\end{displaymath}
which we will express as a linear form in logarithms. For this 
purpose notice that 
\begin{displaymath}
\be = \pm \mu \epsilon_{1}^{a_{1}} \dotsm \epsilon_{d-1}^{a_{d-1}} 
\end{displaymath}
where $a_{1}, \dots, a_{d-1} \in \Z$ and $\mu \in {\mathcal M}$. 
Using the notation,  
\begin{eqnarray*}
\al_{0} = \left| 
	     \frac{\xi^{(i)}-\xi^{(j)}}{\xi^{(i)}-\xi^{(k)}}
	     \frac{\mu^{(k)}}{\mu^{(j)}} \right| 
\text{ and } 
\al_{i} = \left| \frac{\epsilon_{i}^{(k)}} 
			 {\epsilon_{i}^{(j)}} \right|, 
\text{ for $1 \leq i \leq d-1$, }
\end{eqnarray*}
we have 
\begin{displaymath}
\Lambda(i,j,k,\mu) = \log \al_{0} + a_{1} \log \al_{1}      
		     + \dots + a_{d-1} \log \al_{d-1}.  
\end{displaymath}
Let $A = \max \left( 6, |a_{1}|, \dots, |a_{d-1}| \right)$.  
The use of $6$ here is somewhat arbitrary. To apply the 
known results to determine a lower bound for the absolute 
value of these linear forms, we need only have $3$ here. 
However, using $6$ will yield some stronger and simpler 
estimates in what follows, while imposing no significant 
restrictions. 

\begin{lemma}
\label{lem:start}
There exist effectively computable constants $C_{5},C_{6},Y_{1}$ 
and $Y_{2}'$ depending only on $g(X)$, $\K$ and the elements of 
${\mathcal M}$ such that the following statements are true. 

{\rm (i)} If $|y| > Y_{1}$ then $x/y$ is a convergent from the 
continued fraction expansion of $\xi^{(i_{0})}$ for some 
$1 \leq i_{0} \leq d$. 

{\rm (ii)} Suppose that $|y| > Y_{2}'$. Then, for $i_{0}$ as in  
{\rm (i)}, any choice of $j$ and $k$ satisfying 
$i_{0} \neq j \neq k \neq i_{0}$ and any $\mu \in {\mathcal M}$, we have  
\begin{displaymath}
\left| \Lambda(i_{0},j,k,\mu) \right| 
< C_{6} \exp \left( - \frac{d}{C_{5}} A \right).   
\end{displaymath}
\end{lemma}

\begin{proof}
This is a combination of Lemmas 1.1, 1.2, 2.1 and 2.2 of \cite{TW}.  
\end{proof}

Let $\gamma_{1}, \dots, \gamma_{r}$ be algebraic numbers 
with $D=[\Q(\gamma_{1}, \dots, \gamma_{r}):\Q]$ and 
\begin{displaymath}
L = b_{1} \log \gamma_{1} + \dots + b_{r} \log \gamma_{r}, 
\end{displaymath}
with $b_{1}, \dots, b_{r} \in \Z$ and 
$B = \max ( |b_{1}|, \dots, |b_{r}|, 3)$. We put  
\begin{displaymath}
h'(\gamma_{i}) 
= \max \left( h(\gamma_{i}), \left| \log \gamma_{i} \right| / D,  
	      1/D \right) 
\end{displaymath}
for $i=1, \dots, r$ where $h(\gamma_{i})$ is the absolute
logarithmic height of $\gamma_{i}$ and let 
$H = h'(\gamma_{1}) \cdot \dotsm \cdot h'(\gamma_{r})$. 
With $K_{4} = 18(r+1)! r^{r+1} (32D)^{r+2} \log (2Dr) H$, we can 
now state a recent result of Baker \& W\"{u}stholz. 

\begin{lemma}
\label{lem:lform}
If $L \neq 0$ then 
\begin{displaymath}
|L| > \exp \left( -K_{4} \log B \right). 
\end{displaymath}
\end{lemma}

\begin{proof}
See Baker and W\"{u}stholz \cite{BW}. 
\end{proof}

One might expect that we will apply this result to the  
$\Lambda(i_{0},j,k,\mu)$'s. However, it will frequently be 
the case that there are multiplicative relations between 
the $\al_{i}$'s. These relations will allow us to eliminate 
some of the terms in the $\Lambda(i_{0},j,k,\mu)$'s, obtaining 
new linear forms, $\Lambda'(i_{0},j,k,\mu)$, in fewer terms. 
It will turn out that $\Lambda'(i_{0},j,k,\mu)
= t_{0} \Lambda(i_{0},j,k,\mu)$, where $t_{0}$ is a positive 
integer which depends on the particular form. Also the maximum,  
$A'$, of the absolute values of the coefficients of the 
$\Lambda'(i_{0},j,k,\mu)$ can be shown to be at most $A^{2}$ since 
$A \geq 6$ (see Section 4.3). 

Tzanakis and de Weger show in the proof of Lemma 2.4 of 
their paper \cite{TW} that, for $|y| > Y_{2}'$, 
$\Lambda(i_{0},j,k,\mu) \neq 0$. Therefore, we can apply 
Lemma~\ref{lem:lform} to show that, for $|y| > Y_{2}'$, 
\begin{displaymath}
|\Lambda'(i_{0},j,k,\mu)| > \exp \left( -K_{4} \log A' \right).   
\end{displaymath}

Thus, for $|y| > Y_{2}'$,  
\begin{displaymath}
|\Lambda(i_{0},j,k,\mu)| > \exp \left( -2K_{4} \log A \right)/t_{0}.     
\end{displaymath}

In accordance with the notation in \cite{TW}, we shall label 
this constant $2K_{4}$ as $C_{7}$. 

We also mention that the maximal real subfield of a cyclotomic 
field is a Galois extension of $\Q$, so $D=d$, with the exception 
of the equation arising from $n=12,k=-2$ where we let $D=24$. 

We can use this lower bound to prove the following lemma. 

\begin{lemma}
\label{lem:uba}
Put 
\begin{displaymath}
C_{9} = \frac{2C_{5}}{d} \left( \log (t_{0}C_{6}) 
	+ C_{7} \log \left( \frac{C_{5}C_{7}}{d} \right) \right). 
\end{displaymath}

If $|y| > Y_{2}'$ then $A < C_{9}$. 
\end{lemma}

\begin{proof}
This is Lemma 2.4 of \cite{TW} except that we replace $C_{6}$ by 
$t_{0}C_{6}$.  
\end{proof}

Applying this lemma to our linear forms, we obtain an upper bound 
for $A$. By examining the tables at the end of the paper, we see 
that this bound is very large, too large in fact to allow us to 
completely determine all solutions of $F(X,Y)=m$ in any naive way. 
We will use the so-called $L^{3}$ algorithm for lattice basis 
reduction to reduce the size of these upper bounds to the point 
where a direct search is feasible. Our implementation of this 
algorithm will be the iterative integer version given by de Weger 
in Section 3 of \cite{deW1}. We then use either Proposition 3.1 or 
Proposition 3.2 of \cite{TW} to obtain an improved upper bound for 
$A$.  

\subsection{An Improved Upper Bound from the $L^{3}$ Algorithm}

Before using the $L^{3}$ algorithm, there are two conditions we 
must check. First we check for linear relations over $\Z$ among 
$\log \al_{1}, \dots, \log \al_{d-1}$. If such relations exist 
then $|{\bold{b}}_{1}|$ defined below will likely be too small 
for the hypotheses of Lemmas~\ref{lem:l3nz} and \ref{lem:l3z} to 
hold. Therefore we must first eliminate such dependencies. This 
has an advantage too, as we can replace $\Lambda(i_{0},j,k,\mu)$ 
by a linear form with fewer terms and hence obtain smaller values 
of $C_{7}$ and $C_{9}$. Also the $L^{3}$ algorithm will run faster 
since we can then apply it to a smaller matrix with smaller entries. 

We must also check whether $\log \al_{0}$ is a linear combination 
of $\log \al_{1}, \dots, \log \al_{d-1}$ over $\Z$. If this is the 
case, then the quantity $\| s_{k} \|$ in Lemma~\ref{lem:l3nz} will 
be too small. We discuss this  situation after Lemma~\ref{lem:l3nz}. 

By reordering the $\al_{i}$'s if necessary, we may assume 
that $\{ \log \al_{1}, \dots, \log \al_{p} \}$, where 
$p \leq d-1$, is a $\Q$--linearly independent set. Proceeding as 
Tzanakis and de Weger do in case (iii) of Section II.3 of their 
paper \cite{TW}, we can find integers $t_{0} > 0$ and $t_{ij}$ 
for $1 \leq i \leq p, p+1 \leq j \leq d-1$ such that 
\begin{displaymath}
t_{0} \log \al_{j} = \sum_{i=1}^{p} t_{ij} \log \al_{i} 
\text{ for $j=p+1, \dots, d-1$} 
\end{displaymath}
(notice that we use $t_{0}$ and $t_{ij}$ where Tzanakis and 
de Weger \cite{TW} use $d$ and $d_{ij}$). Using these relations, 
we can eliminate the terms for $\al_{p+1}, \dots, \al_{d-1}$ from 
our linear form $\Lambda(i_{0},j,k,\mu)$, obtaining 
\begin{displaymath}
\Lambda'(i_{0},j,k,\mu) = t_{0} \Lambda(i_{0},j,k,\mu) 
= t_{0} \log \al_{0} + \sum_{i=1}^{p} a_{i}' \log \al_{i},  
\end{displaymath}
where 
\begin{displaymath}
a_{i}' = t_{0} a_{i} + \sum_{j=p+1}^{d-1} t_{ij} a_{j}.  
\end{displaymath}

Letting $T=\max(t_{0},|t_{ij}|:1 \leq i \leq p, p+1 \leq j \leq d-1)$, 
we have 
\begin{equation}
\label{eq:upbnd}
|\Lambda'(i_{0},j,k,\mu)| < t_{0} C_{6} \exp \left( \frac{d}{C_{4}} A \right),  
\text{ with } A < C_{9}  
\end{equation}
by Lemmas~\ref{lem:start}(ii)--\ref{lem:uba}, and 
\begin{displaymath}
|a_{i}'| \leq (d-p)TA \text{ for $i=1, \dots, p$}.  
\end{displaymath}

With $\log \al_{1}, \dots, \log \al_{p}$ being $\Q$--linearly 
independent we now apply the $L^{3}$ algorithm to the matrix 
\begin{displaymath}
{\mathcal A} = \left( \begin{matrix}
1                           & 0      & \dotsc & 0 \\
0                           & \ddots & 0      & 0 \\
\vdots                      & 0      & 1      & 0 \\
{[c_{0} \log \al_{1}]}      & \dotsc & \dotsc & [c_{0} \log \al_{p}] 
		  \end{matrix} 
	   \right) 
\end{displaymath}
where $c_{0}$ is a real number somewhat larger than 
$C_{9}^{p}$ and obtain a matrix we will denote by ${\mathcal B}$. 
Let 
\begin{displaymath}
\bold{x} = (0, \dots, 0, -[c_{0} t_{0} \log \al_{0}])^{T} 
	 = \sum_{i=1}^{p} s_{i} {\bold{b}}_{i} 
\end{displaymath}
where ${\bold{b}}_{i}$ is the vector formed from the $i$-th 
column of ${\mathcal B}$ and $s_{1}, \dots, s_{p} \in \R$. Let 
$k$ be the largest integer such $s_{k} \not\in \Z$ and for 
$x \in \R$, denote by $\| x \|$ the distance from $x$ to the 
nearest integer. 

\begin{lemma}
\label{lem:l3nz}
Suppose 
\begin{displaymath}
2^{-(p-1)/2} \| s_{k} \| \cdot \left| {\bold{b}}_{1} \right|  
\geq \sqrt{p^{2}+5p+3} \, (d-p)TC_{9}. 
\end{displaymath}
Then there are no solutions of 
$|\Lambda(i_{0},j,k,\mu)| < C_{6} \exp(-dA/C_{5})$ with 
\begin{displaymath}
A > \frac{C_{5}}{d} 
    \log \left( \frac{c_{0}t_{0}C_{6}}{(d-p)TC_{9}} \right). 
\end{displaymath}
\end{lemma}
					       
\begin{proof}
This lemma is a slight modification of Prop.\ 3.2 of \cite{TW} 
(we obtain a result for $A$ whereas their result pertains to 
an upper bound for the $|a_{i}'|$'s). The proof is nearly 
identical to the proof of Lemma 3.10 in \cite{deW2}, the 
result upon which Prop.\ 3.2 of \cite{TW} is based, except 
that at the very end of the proof we use the upper bound for 
$|\Lambda'(i_{0},j,k,\mu)|$ in (\ref{eq:upbnd}) which is in 
terms of $A$ to get our result in terms of $A$. 
\end{proof}

As we stated above, if $\log \al_{0}$ is a linear combination 
of $\log \al_{1}, \dots, \log \al_{d-1}$ over $\Z$ then the 
quantity $\| s_{k} \|$ may be quite small. In this case we 
eliminate $\log \al_{0}$ from our linear forms. This situation 
only occurs for the Thue equations we consider with $n=7$ 
and $n=9$.  Let us now describe how we deal with this situation.  

Suppose 
\begin{displaymath}
t_{0}\log\al_{0} + \dots + t_{d-1}\log\al_{d-1} = 0, 
\end{displaymath}
with $t_{0}, \dots, t_{d-1} \in \Z$ satisfying 
$t_{0} > 0$, and let $T = \max (|t_{0}|, \dots, |t_{d-1}|)$. 
Letting $a_{i}' = a_{i}t_{0}-t_{i}$ for $1 \leq i \leq d-1$, 
we have 
\begin{displaymath}
\Lambda'(i_{0},j,k,\mu) = t_{0} \Lambda(i_{0},j,k,\mu) 
=  a_{1}' \log \al_{1} + \dots + a_{d-1}' \log \al_{d-1}.   
\end{displaymath}
Notice that $A' = \max ( 3, |a_{1}'|, \dots, |a_{d-1}'| ) 
\leq T(A+1) < 1.17TA$, since $A \geq 6$. 

Thus, we have 
\begin{displaymath}
|\Lambda'(i_{0},j,k,\mu)| 
< t_{0}C_{6} \exp \left( - \frac{6d}{7C_{5}T} A' \right),    
\text{ with } A' < 1.17TC_{9},     
\end{displaymath}
by Lemmas~\ref{lem:start}(ii)--\ref{lem:uba}. 

For the Thue equations arising from $n=7$ and $9$, 
$\log \al_{1}, \dots, \log \al_{d-1}$ are linearly independent 
over $\Q$. Thus we are ready to apply the $L^{3}$ algorithm, 
our use of which shall be similar to the previous case: we 
apply the algorithm to the same matrix ${\mathcal A}$ as before, 
here with $p=d-1$, obtaining the matrix ${\mathcal B}$. 

\begin{lemma}
\label{lem:l3z}
Suppose 
\begin{displaymath}
2^{(d-2)/2} \left| {\bold{b}}_{1} \right| 
> 1.17 \sqrt{d^{2}+d-2} \, TC_{9}. 
\end{displaymath}
Then there are no solutions of 
$|\Lambda(i_{0},j,k,\mu)| < C_{6} \exp (-dA/C_{5})$ 
with   
\begin{displaymath}
A > \left( \frac{1.16C_{5}T}{d} 
    \log \left( \frac{0.85c_{0}C_{6}}{C_{9}} \right) + T \right) / t_{0}. 
\end{displaymath}
\end{lemma}
					       
\begin{proof}
Letting $q=d-1$, we have $q^{2}+q-1=d^{2}-d-1$ and so the 
hypothesis of Proposition 3.1 from \cite{TW} holds. Thus,  
we see that 
\begin{displaymath}
A' > \frac{1.16C_{5}T}{d} 
     \log \left( \frac{c_{0}t_{0}C_{6}}{1.17TC_{9}} \right). 
\end{displaymath}

Noting that $T \geq t_{0}$, our bound for $A$ now follows 
from the fact that our definition of the $a_{i}'$'s shows 
that $A \leq (A'+T)/t_{0}$. 
\end{proof}

\subsection{Searching for Solutions of $F(X,Y)=m$}

We construct linear forms in logarithms as in Section 4.1 
and using the methods described above get a good upper bound 
for the size of the coefficients of these linear forms. We 
must do this for each $i_{0}$ between $1$ and $d$ and each 
$\mu \in {\mathcal M}$. However, because of the large bound we 
obtain from Lemma~\ref{lem:uba}, the $L^{3}$ algorithm takes 
a long time to run when $d$ is large. But one notices that 
only $\al_{0}$ depends on $i_{0}$ and $\mu$ so the matrix 
${\mathcal A}$ defined in Section 4.2 depends only on the $j$ 
and $k$ defined in Section 4.1. When $d \geq 4$, we can 
choose $(j_{1},k_{1})$ and $(j_{2},k_{2})$ with 
$(j_{1},k_{1}) \neq (j_{2},k_{2})$ so that we need only apply 
the $L^{3}$ algorithm to two matrices of the form ${\mathcal A}$: 
if $i_{0}=j_{1}$ or $i_{0}=k_{1}$ we let $j=j_{2}$ and $k=k_{2}$ 
otherwise we let $j=j_{1}$ and $k=k_{1}$. We choose 
$(j_{1},k_{1})$ and $(j_{2},k_{2})$ so that 
$\epsilon_{i}^{(k_{1})}/\epsilon_{i}^{(j_{1})}$ and 
$\epsilon_{i}^{(k_{2})}/\epsilon_{i}^{(j_{2})}$ are conjugates 
for $i=1, \dots, d-1$. This simplifies the height calculations 
necessary for Lemma~\ref{lem:lform}. And more importantly, we 
choose them to minimize the number of $\Q$--linearly independent 
$\log \al_{i}$'s. 

We then compute the quantities $C_{5},C_{6},Y_{1}$ and $Y_{2}'$.   
Next we apply Lemma~\ref{lem:lform} to the $\Lambda'(i_{0},j,k,\mu)$'s 
to determine $K_{4}$. Notice that this lemma will give us a lower 
bound in terms of $A'$. We saw that in Section 4.2 that either 
$A' \leq (d-p)TA$ or $A' \leq 1.1TA$. An examination of the relations  
in the next section shows that both $(d-p)T$ and $1.1T$ are at most 
6 and, since we have assumed $A \geq 6$, we have $A' \leq A^{2}$. 
Thus for $C_{7}$ we use $2K_{4}$. From these quantities we find 
$C_{9}$. Applying the $L^{3}$ algorithm as described above 
considerably reduces the upper bound for $A$ and then applying 
the $L^{3}$ a second time using this new upper bound for $A$ in 
place of $C_{9}$, we obtain a still smaller upper bound. At this 
point we wish to determine an upper bound for $|y|$ from this 
last upper bound for $A$. We use the following lemma. 

\begin{lemma}
\label{lem:uby}
Suppose $(x,y) \in {\Z}^{2}$ is a solution of $F(X,Y)=m$ and $A < C_{10}$. 
Then 
\begin{displaymath}
|y| \leq Y_{3} = \min_{1 \leq j_{1} < j_{2} \leq d} 
\left( \mu_{+} \frac{E_{j_{1}}^{C_{10}}+E_{j_{2}}^{C_{10}}} 
		    {\xi^{(j_{1})} - \xi^{(j_{2})}} \right)
\end{displaymath}
where $\displaystyle E_{j} = \prod_{i=1}^{d-1} 
{\left| \epsilon_{i}^{(j)} \right|}^{v_{ij}}$,  
$\displaystyle 
\mu_{+} = \max_{1 \leq i \leq d, \mu \in {\mathcal M}} 
	  \left| \mu^{(i)} \right|$ 
and $v_{ij} = \pm 1$ whichever one makes 
${ \left| \epsilon_{i}^{(j)} \right| }^{v_{ij}} \geq 1$. 
\end{lemma}

\begin{proof}
This is proven on page 118 of the paper \cite{TW} of Tzanakis 
and de Weger. 
\end{proof}

Now we perform a direct search for solutions with 
$|y| \leq Y_{1}$ and then check whether $(x,y)$ is a solution 
of $F(X,Y)=m$ where $x/y$ is a convergent of $\xi^{(i)}$ with 
$|y| \leq Y_{3}$ for each $1 \leq i \leq d$. In this manner we 
are able to determine the complete solution of the Thue equation 
$F(X,Y)=m$. 

\section{Dependence Relations}

In Section 4.2, we described what to do when dependence 
relations arise among the $\log \al_{i}$'s. Here we give 
the relations that were found in our applications. These 
relations were found either by making use of the nice form 
of the $\al_{i}$'s when $n$ is a prime power and by direct 
search otherwise. 

Let us first establish an ordering of the units and their 
conjugates. In the case of $n \neq 12$, let 
$a_{1}=1, \dots, a_{\varphi(n)/2-1}$ be the increasing sequence 
of positive integers less than $n/2$ which are relatively prime 
to $n$. In accordance with Lemma~\ref{lem:units}, we let 
\begin{displaymath}
\epsilon_{i}^{(j)} = \pm \frac{\sin(a_{i+1}a_{j} \pi /n)}
			  {\sin(a_{j} \pi /n)}
\end{displaymath}
denote the $j$-th conjugate of $i$-th fundamental unit. Notice  
that for our purposes here, knowledge of the conjugates up to 
sign suffices, for the $\alpha_{i}$'s are defined to be the 
absolute value of quotients of these conjugates. 

For the Thue equation which arises from $n=12$, we order the 
roots of $F(S,1)$ as follows: 
$\alpha^{(1)} = 1+\sqrt{3}+\sqrt{6+2\sqrt{3}}, 
\alpha^{(2)}  = 1+\sqrt{3}-\sqrt{6+2\sqrt{3}}, 
\alpha^{(3)}  = 1-\sqrt{3}+\sqrt{6-2\sqrt{3}}$ and  
$\alpha^{(4)} = 1-\sqrt{3}-\sqrt{6-2\sqrt{3}}$.  
The meaning of $\epsilon_{i}^{(j)}$ is then clear using the 
labeling of the fundamental units in Section 3.4. 

We now consider the dependence relations themselves. 

For $n=7$ and $9$, we let $(j,k)=(2,3)$ when $i_{0}=1$, 
$(j,k)=(1,3)$ when $i_{0}=2$, and $(j,k)=(1,2)$ when 
$i_{0}=3$. As we mentioned in Section 4.2, 
$\log \al_{0},\log \al_{1}$ and $\log \al_{2}$ are 
linearly dependent over $\Q$. When $n=7$ and $m=\pm 1$, 
we have 
\begin{displaymath}
3 \log \al_{0} = \log \al_{1} - 2 \log \al_{2}, 
\end{displaymath}
for each choice of $i_{0},j$ and $k$. When $n=7$ and 
$m=\pm 7$, we have 
\begin{displaymath}
\log \al_{0} = -\log \al_{1} - 2 \log \al_{2}, 
\end{displaymath}
for each choice of $i_{0},j$ and $k$. So we let $t_{0}=3$ 
and $T=6$ 

When $n=9$ and $m=\pm 1$, we have 
\begin{displaymath}
3 \log \al_{0} = 2 \log \al_{1} - \log \al_{2}, 
\end{displaymath}
for each choice of $i_{0},j$ and $k$. When $n=9$ and 
$m=\pm 3$, we have 
\begin{displaymath}
\log \al_{0} = -\log \al_{2}, 
\end{displaymath}
for each choice of $i_{0},j$ and $k$. Here we can use 
$t_{0}=3$ and $T=3$.

For $n=11$, we let $(j,k)=(1,2)$ if $i_{0} \neq 1,2$ and  
$(j,k)=(3,5)$ otherwise. With this choice, 
$\log \al_{1}, \log \al_{2}, \log \al_{3}, \log \al_{4}$ 
are $\Q$--linearly independent.  

For $n=12$ with $k=-2$, we let $(j,k) = (1,3)$ if 
$i_{0} \neq 1,3$ and $(j,k) = (2,4)$ otherwise. In 
both cases, $\log \al_{1}, \log \al_{2}, \log \al_{3}$ 
are $\Q$--linearly independent, and so we use 
let $\Lambda' = \Lambda$.  

For $n=13$, we let $(j,k) = (1,5)$ if $i_{0} \neq 1,5$ and 
$(j,k) = (2,3)$ otherwise. In both cases, we have 
\begin{displaymath}
\al_{1} = \al_{3} \al_{5} / \al_{2} 
\hspace{3 mm} \text{ and } \hspace{3 mm} 
\al_{4} = \al_{3} \al_{5}. 
\end{displaymath}

We let 
\begin{displaymath}
\Lambda' = \log \al_{0} + a_{1}' \log \al_{2}
	   + a_{2}' \log \al_{3} + a_{3}' \log \al_{5}.  
\end{displaymath}

For $n=15$, we let $(j,k)=(1,2)$ if $i_{0} \neq 1,2$ and  
$(j,k)=(3,4)$ otherwise. With this choice, 
$\log \al_{1}, \log \al_{2}, \log \al_{3}$ are 
$\Q$--linearly independent.  

For $n=16$, we let $(j,k) = (1,4)$ if $i_{0} \neq 1,4$ and 
$(j,k) = (2,3)$ otherwise. In both cases, we have 
$\al_{3} = \al_{1} \al_{2}$. So we let 
\begin{displaymath}
\Lambda' = \log \al_{0} + a_{1}' \log \al_{1} + a_{2}' \log \al_{2}. 
\end{displaymath}

For $n=17$, we let $(j,k) = (1,4)$ if $i_{0} \neq 1,4$ and 
$(j,k) = (2,8)$ otherwise. In both cases, we have 
\begin{displaymath} 
\al_{1} = \al_{5} \al_{6} / \al_{7}, 
\hspace{3 mm} 
\al_{2} = \al_{5} \al_{6} / \al_{4} 
\hspace{3 mm} \text{ and } \hspace{3 mm} 
\al_{3} = \al_{5} \al_{6}. 
\end{displaymath}

So we let 
\begin{displaymath}
\Lambda' = \log \al_{0} + a_{1}' \log \al_{4} 
	   + a_{2}' \log \al_{5} + a_{3}' \log \al_{6} 
	   + a_{4}' \log \al_{7}.  
\end{displaymath}

For $n=19$, we let $(j,k) = (1,7)$ if $i_{0} \neq 1,7$ and 
$(j,k) = (2,5)$ otherwise. In both cases, we have 
\begin{displaymath} 
\al_{3} = \al_{1} \al_{2} \al_{4} 
		 / \left( \al_{5} \al_{8} \right) 
\hspace{3 mm}  \text{ and }  \hspace{3 mm} 
\al_{6} = \al_{1} \al_{2} \al_{4} / \al_{7}.  
\end{displaymath}

So we let 
\begin{displaymath}
\Lambda' = \log \al_{0} + a_{1}' \log \al_{1}
	   + a_{2}' \log \al_{2} + a_{3}' \log \al_{4}
	   + a_{4}' \log \al_{5} + a_{5}' \log \al_{7} 
	   + a_{6}' \log \al_{8}.  
\end{displaymath}

For $n=20$, we let $(j,k)=(1,2)$ if $i_{0} \neq 1,2$ and  
$(j,k)=(3,4)$ otherwise. With this choice, 
$\log \al_{1},\log \al_{2}$ and $\log \al_{3}$ are $\Q$-linearly 
independent, so we let $\Lambda'=\Lambda$. 

For $n=21$, we let $(j,k)=(1,5)$ if $i_{0} \neq 1,5$ and 
$(j,k)=(2,4)$ otherwise. In both cases, we have 
\begin{displaymath}
\al_{4} = \al_{3}/\al_{1} 
\hspace{3 mm} \text{ and } \hspace{3 mm} 
\al_{5} = 1 / \al_{1}.  
\end{displaymath}

So we let  
\begin{displaymath}
\Lambda' = \log \al_{0} + a_{1}' \log \al_{1}
	   + a_{2}' \log \al_{2} + a_{3}' \log \al_{3}.  
\end{displaymath}

For $n=23$, we let $(j,k)=(1,2)$ if $i_{0} \neq 1,2$ and  
$(j,k)=(3,6)$ otherwise. With this choice, 
$\log \al_{1}, \dots, \log \al_{10}$ are $\Q$--linearly independent.  

For $n=24$, we let $(j,k) = (1,2)$ if $i_{0} \neq 1,2$ and 
$(j,k) = (3,4)$ otherwise. In both cases, we have 
$\al_{2} \al_{3} = \al_{1}^{3}$. So   
\begin{displaymath}
\Lambda' = 3 \log \al_{0} + a_{2}' \log \al_{2} + a_{3}' \log \al_{3}. 
\end{displaymath}

For $n=25$, we let $(j,k) = (2,9)$ if $i_{0} \neq 2,9$ and 
$(j,k) = (4,3)$ otherwise. In both cases, we have 
\begin{eqnarray*} 
\al_{1} = \al_{5} / \al_{8}, &              & \al_{2} = \al_{5} / \al_{3}, \\ 
\al_{4} = \al_{5} / \al_{6}  & \text{ and } & \al_{7} = \al_{5} / \al_{9}. 
\end{eqnarray*}

So we let 
\begin{displaymath}
\Lambda' = \log \al_{0} + a_{1}' \log \al_{3}
	   + a_{2}' \log \al_{5} + a_{3}' \log \al_{6} 
	   + a_{4}' \log \al_{8} + a_{5}' \log \al_{9}.  
\end{displaymath}

For $n=29$, we let $(j,k) = (1,12)$ if $i_{0} \neq 1,12$ and 
$(j,k) = (2,5)$ otherwise. In both cases, we have 
\begin{eqnarray*} 
\al_{1}  = \al_{10} \al_{12} / \al_{4}, & & 
\al_{2}  = \al_{10} \al_{12} / \al_{6}, \\ 
\al_{3}  = \al_{10} \al_{12} / \al_{9}, & &  
\al_{5}  = \al_{10} \al_{12} / \al_{13}, \\ 
\al_{7}  = \al_{10} \al_{12} / \al_{8} & \text{ and } & 
\al_{11} = \al_{10} \al_{12}. 
\end{eqnarray*}

So we let  
\begin{eqnarray*}
\Lambda' & = & \log \al_{0} + a_{1}' \log \al_{4}
	       + a_{2}' \log \al_{6} + a_{3}' \log \al_{8} \\ 
	 & &   + a_{4}' \log \al_{9} + a_{5}' \log \al_{10} 
	       + a_{6}' \log \al_{12} + a_{7}' \log \al_{13}.  
\end{eqnarray*}

\section{Tables of Results}

\subsection{Equations Solved by the Algorithm of Tzanakis \& de Weger}

For $n=7,9,11$,$13$,$15$,$16$,$17$,$19$,$20$,$21$,$23$,$24$,
$25$ and $29$, as well as $n=12$ with $k=-2$, we used the method 
of Tzanakis and de Weger as described in Section 4 to solve the 
Thue equations which arise. The method was implemented using the 
MAPLE V Computer Algebra System on an 80486 DX2 based 
IBM--compatible PC running at 50 MHz. In the first two tables, 
we list the equations solved by this method. These are followed 
by three tables containing an abridgement of the output from 
these programs for each equation solved. The entries in these 
latter tables have been rounded up or down, as appropriate. 
Requests for more information regarding these computations 
are, of course, welcome. 

The first entries in Tables 5--7, before $d_{1}$, are listed 
using the notation of Section 4. As mentioned above, the $L^{3}$ 
algorithm was used twice. The first time we let $c_{0} = 10^{d_{1}}$ 
from which we obtained $A \leq A_{1}$. The second time, we let 
$c_{0} = 10^{d_{2}}$ and found that $A \leq A_{2}$. These quantities, 
$A_{1},A_{2},d_{1}$ and $d_{2}$, are listed in these tables. Finally, 
we used Lemma~\ref{lem:uby} to obtain an upper bound for $|y|$ from 
$A \leq A_{2}$. This is listed in Tables 5--7 under the entry $Y_{3}$. 
$X_{4}$ (resp. $Y_{4}$) is the maximum of the absolute value of $x$ 
(resp. $y$) for all solutions $(x,y)$ to the Thue equations which 
arise for each $n$. A complete list of solutions has been omitted 
to save space, however, $X_{4}$ and $Y_{4}$ are sufficiently small 
that the interested reader could easily determine all solutions. The 
last entry gives the CPU time used for each $n$. We believe that the 
time is of interest showing as it does that the method of Tzanakis 
and de Weger is practical even for Thue equations of moderate degree, 
provided the necessary system of fundamental units and factorization 
of $m$ are known. 

\begin{table}
\label{tab:eqns1}
\caption{}
\begin{center}
\begin{tabular}{||l|l||}                               
\hline
$n$ &  $F_{n}(X,Y) = m$                             
\\ \hline 
7   &  $X^{3}+X^{2}Y-2XY^{2}-Y^{3} = \pm 1, \pm 7$  
\\ \hline
9   &  $X^{3}-3XY^{2}+Y^{3} = \pm 1, \pm 3$         
\\ \hline
11  &  $X^{5}+X^{4}Y-4X^{3}Y^{2}-3X^{2}Y^{3}+3XY^{4}+Y^{5} = \pm 1, \pm 11$         
\\ \hline
12  &  $X^{4}-4X^{3}Y-12X^{2}Y^{2}+8XY^{3}+4Y^{4} = 1$         
\\ \hline
13  &  $X^{6}+X^{5}Y-5X^{4}Y^{2}-4X^{3}Y^{3}+6X^{2}Y^{4}+3XY^{5}-Y^{6} 
= \pm 1, \pm 13$         
\\ \hline
15  &  $X^{4}-X^{3}Y-4X^{2}Y^{2}+4XY^{3}+Y^{4} = \pm 1, \pm 5$         
\\ \hline
16  &  $X^{4}-4X^{2}Y^{2}+2Y^{4} = \pm 1, \pm 2$         
\\ \hline
17  &  $X^{8}+X^{7}Y-7X^{6}Y^{2}-6X^{5}Y^{3}+15X^{4}Y^{4}$ \\ \cline{2-2} 
    &  $+10X^{3}Y^{5}-10X^{2}Y^{6}-4XY^{7}+Y^{8} 
= \pm 1, \pm 17$  
\\ \hline
\end{tabular}
\end{center}
\end{table}

\begin{table}
\begin{center}
\label{tab:eqns2}
\caption{}
\begin{tabular}{||l|l||}                               
\hline
$n$ &  $F_{n}(X,Y) = m$                             
\\ \hline 
19  &  $X^{9}+X^{8}Y-8X^{7}Y^{2}-7X^{6}Y^{3}+21X^{5}Y^{4}$ \\ \cline{2-2} 
    &  $+15X^{4}Y^{5}-20X^{3}Y^{6}-10X^{2}Y^{7}+5XY^{8}+Y^{9} 
= \pm 1, \pm 19$         
\\ \hline
20  &  $X^{4}-5X^{2}Y^{2}+5Y^{4} = \pm 1, \pm 5$         
\\ \hline
21  &  $X^{6}-X^{5}Y-6X^{4}Y^{2}+6X^{3}Y^{3}+8X^{2}Y^{4}-8XY^{5}+Y^{6} 
= \pm 1, \pm 7$         
\\ \hline
23  &  $X^{11}+X^{10}Y-10X^{9}Y^{2}-9X^{8}Y^{3}+36X^{7}Y^{4}+28X^{6}Y^{5}$ 
       \\ \cline{2-2}  
    &  $-56X^{5}Y^{6}-35X^{4}Y^{7}+35X^{3}Y^{8}+15X^{2}Y^{9}-65XY^{10}-Y^{11} 
= \pm 1, \pm 23$         
\\ \hline
24  &  $X^{4}-4X^{2}Y^{2}+Y^{4} = \pm 1, \pm 2$         
\\ \hline
25  &  $X^{10}-10X^{8}Y^{2}+35X^{6}Y^{4}+X^{5}Y^{5}-50X^{4}Y^{6}$ 
       \\ \cline{2-2}  
    &  $-5X^{3}Y^{7}+25X^{2}Y^{8}+5XY^{9}-Y^{10} = \pm 1, \pm 5$         
\\ \hline
29  &  $X^{14}+X^{13}Y-13X^{12}Y^{2}-12X^{11}Y^{3}+66X^{10}Y^{4}$ 
       \\ \cline{2-2}  
    &  $+55X^{9}Y^{5}-165X^{8}Y^{6}-120X^{7}Y^{7}+210X^{6}Y^{8}+126X^{5}Y^{9}$ 
    \\ \cline{2-2} 
    &  $-126X^{4}Y^{10}-56X^{3}Y^{11}+28X^{2}Y^{12}+7XY^{13}-Y^{14} 
= \pm 1, \pm 29$ 
\\ \hline         
\end{tabular}
\end{center}
\end{table}

\begin{table}
\begin{center}
\label{tab:res1}
\caption{}
\begin{tabular}{||l|r|r|r|r|r||}                               \hline
$n$             & $7$               & $9$                 & $11$          
		& $12 (k=-2)$       & $13$                  \\ \hline
$d$             & $3$               & $3$                 & $5$           
		& $4$               & $6$                   \\ \hline 
$Y_{1}$         & $49$              & $19$                & $6$           
		& $2$               & $5$                   \\ \hline 
$Y_{2}'$        & $49$              & $19$                & $12$          
		& $2$               & $17$                  \\ \hline
$d/C_{5}$       & $1.127$           & $1.508$             & $1.584$
		& $2.171$           & $1.834$               \\ \hline
$C_{6}$         & $5900$            & $18000$             & $5 \cdot 10^{8}$ 
		& $48000$           & $2 \cdot 10^{11}$     \\ \hline
$H$             & $0.218$           & $0.317$             & $0.0373$ 
		& $1.162$           & $0.132$               \\ \hline   
$C_{7}$         & $8 \cdot 10^{10}$ & $1.2 \cdot 10^{11}$ & $1.6 \cdot 10^{23}$ 
		& $6 \cdot 10^{24}$ & $1.2 \cdot 10^{20}$   \\ \hline 
$C_{9}$         & $4 \cdot 10^{12}$ & $4 \cdot 10^{12}$   & $2 \cdot 10^{25}$ 
		& $3 \cdot 10^{27}$ & $6 \cdot 10^{21}$     \\ \hline
$d_{1}$         & $29$              & $29$                & $120$         
		& $90$              & $77$                  \\ \hline 
$A_{1}$         & $61$              & $84$                & $150$         
		& $72$              & $82$                  \\ \hline 
$d_{2}$         & $8$               & $8$                 & $20$          
		& $18$              & $16$                  \\ \hline
$A_{2}$         & $30$              & $42$                & $38$         
		& $22$              & $31$                  \\ \hline 
$Y_{3}$         & $8 \cdot 10^{13}$ & $3 \cdot 10^{27}$   & $10^{33}$ 
		& $3 \cdot 10^{36}$ & $3 \cdot 10^{32}$     \\ \hline 
$(X_{4},Y_{4})$ & $(9,9)$           & $(3,3)$             & $(2,1)$           
		& $(1,0)$           & $(3,2)$               \\ \hline 
time            & $70$ s            & $62$ s              & $94$ s       
		& $81$ s            & $60$ s                \\ \hline 
\end{tabular}
\end{center}
\end{table}

\begin{table}
\label{tab:res2}
\caption{}
\begin{center}
\begin{tabular}{||l|r|r|r|r|r||}                                 \hline
$n$             & $15$                & $16$                & $17$ 
		& $19$                & $20$                  \\ \hline
$d$             & $4$                 & $4$                 & $8$ 
		& $9$                 & $4$                   \\ \hline
$Y_{1}$         & $9$                 & $4$                 & $4$ 
		& $4$                 & $6$                   \\ \hline          
$Y_{2}'$        & $11$                & $4$                 & $29$ 
		& $37$                & $7$                   \\ \hline
$d/C_{5}$       & $1.034$             & $1.738$             & $2.265$  
		& $2.564$             & $0.869$               \\ \hline 
$C_{6}$         & $610000$            & $30000$             & $7 \cdot 10^{16}$   
		& $6 \cdot 10^{19}$   & $5.6 \cdot 10^{7}$    \\ \hline
$H$             & $0.129$             & $0.671$             & $0.0887$    
		& $0.0207$            & $0.333$               \\ \hline
$C_{7}$         & $8.7 \cdot 10^{18}$ & $5.2 \cdot 10^{15}$ & $1.2 \cdot 10^{25}$  
		& $1.2 \cdot 10^{34}$ & $2.3 \cdot 10^{19}$   \\ \hline 
$C_{9}$         & $8 \cdot 10^{20}$   & $3 \cdot 10^{21}$   & $6 \cdot 10^{26}$  
		& $7 \cdot 10^{35}$   & $3 \cdot 10^{21}$     \\ \hline
$d_{1}$         & $70$                & $48$                & $130$  
		& $240$               & $74$                  \\ \hline
$A_{1}$         & $122$               & $46$                & $121$   
		& $200$               & $159$                 \\ \hline
$d_{2}$         & $13$                & $7$                 & $24$  
		& $37$                & $16$                  \\ \hline
$A_{2}$         & $37$                & $12$                & $38$   
		& $48$                & $57$                  \\ \hline
$Y_{3}$         & $9 \cdot 10^{30}$   & $8 \cdot 10^{10}$   & $5 \cdot 10^{56}$  
		& $9 \cdot 10^{83}$   & $2 \cdot 10^{60}$     \\ \hline
$(X_{4},Y_{4})$ & $(4,3)$             & $(2,1)$             & $(2,1)$  
		& $(2,1)$             & $(2,1)$               \\ \hline
time            & $64$ s              & $34$ s              & $139$ s 
		& $583$ s             & $68$ s                \\ \hline
\end{tabular}   
\end{center}
\end{table}

\begin{table}
\label{tab:res3}
\caption{}
\begin{center}
\begin{tabular}{||l|r|r|r|r|r||}                                 \hline
$n$             & $21$                & $23$                & $24$             
		& $25$                & $29$                  \\ \hline
$d$             & $6$                 & $11$                & $4$             
		& $10$                & $14$                  \\ \hline
$Y_{1}$         & $4$                 & $3$                 & $5$             
		& $4$                 & $3$                   \\ \hline
$Y_{2}'$        & $22$                & $54$                & $5$             
		& $63$                & $85$                  \\ \hline
$C_{6}$         & $7 \cdot 10^{9}$    & $6 \cdot 10^{25}$   & $300000$      
		& $4 \cdot 10^{24}$   & $2 \cdot 10^{35}$     \\ \hline 
$d/C_{5}$       & $0.833$             & $2.89$              & $0.758$        
		& $2.07$              & $3.854$               \\ \hline 
$H$             & $0.260$             & $0.00026$           & $0.903$         
		& $0.125$             & $0.0177$              \\ \hline
$C_{7}$         & $2.3 \cdot 10^{20}$ & $10^{53}$           & $7 \cdot 10^{15}$ 
		& $3.4 \cdot 10^{30}$ & $5.5 \cdot 10^{40}$   \\ \hline 
$C_{9}$         & $3 \cdot 10^{22}$   & $8 \cdot 10^{54}$   & $7 \cdot 10^{17}$ 
		& $3 \cdot 10^{32}$   & $3 \cdot 10^{42}$     \\ \hline 
$d_{1}$         & $76$                & $588$               & $45$          
		& $182$               & $342$                 \\ \hline
$A_{1}$         & $173$               & $445$               & $96$             
		& $193$               & $199$                 \\ \hline 
$d_{2}$         & $16$                & $63$                & $11$             
		& $29$                & $42$                  \\ \hline
$A_{2}$         & $63$                & $68$                & $41$          
		& $56$                & $44$                  \\ \hline 
$Y_{3}$         & $2 \cdot 10^{87}$   & $9 \cdot 10^{149}$  & $8 \cdot 10^{51}$    
		& $3 \cdot 10^{123}$  & $3 \cdot 10^{127}$    \\ \hline 
$(X_{4},Y_{4})$ & $(2,1)$             & $(2,1)$             & $(2,2)$       
		& $(2,1)$             & $(2,1)$               \\ \hline
time            & $169$ s             & $11356$ s           & $51$ s        
		& $595$ s             & $2847$ s              \\ \hline
\end{tabular}                                    
\end{center}
\end{table}

\subsection{Solutions of the Equations with $n \equiv 2 \bmod 4$ 
or nonsquarefree} 

Using the notation of Lemma~\ref{lem:sqft}, if $n \geq 4$ is not 
a power of three with $m=3$, then $P(n/(3,n)) = P(m/(3,m))$. If 
$\varphi (m)/2 \geq 3$ and we have determined all solutions of 
$F_{m}(X,Y) = \pm 1, \pm P(m/(3,m))$ then we can use 
Lemma~\ref{lem:sqft}(ii) to find all solutions of 
$F_{n}(X,Y) = \pm 1, \pm P(n/(n,3))$. 

Similarly, if $n=2m$ where $m$ is odd then we use Lemma~\ref{lem:sqft}(i). 
In Table 8 we give the values of $m,X_{4}$ and $Y_{4}$, 
where $X_{4}$ and $Y_{4}$ are as in the previous tables, when 
$14 \leq n \leq 30$ is non-squarefree or $n \equiv 2 \bmod 4$.  

\begin{table}
\label{tab:sqr}
\caption{}
\begin{center}
\begin{tabular}{||c|c|c|c|c|c|c|c||}                       \hline
$n$             &  $14$   &  $18$   &  $22$   &  $26$    
		&  $27$   &  $28$   &  $30$             \\ \hline 
$m$             &   $7$   &   $9$   &  $11$   &  $13$    
		&   $9$   &  $14$   &  $15$             \\ \hline 
$(X_{4},Y_{4})$ & $(9,9)$ & $(3,3)$ & $(2,1)$ & $(3,2)$ 
		& $(2,1)$ & $(2,1)$ & $(4,3)$           \\ \hline 
\end{tabular}
\end{center}
\end{table}

\end{document}